\documentclass{article}
\usepackage{bbm}
\usepackage{amsmath,amssymb,amsfonts}
\author{H\"useyin \c{C}akalli*, \c{C}i\u{g}dem G\"und\"uz Aras**, and Ayse Sonmez***\\
$*$Maltepe University, Marmara E\u{g}itim K\"oy\"u, TR 34857,
\.{I}stanbul, Turkey \\ Phone:$(+90216)6261050$ ext:2248, \,
fax:$(+90216)6261113$, \\ e-mail:hcakalli@maltepe.edu.tr;
hcakalli@gmail.com \\ $**$ Kocaeli University,
Department of Mathematics, Kocaeli, Turkey \\
Phone:$(+90262)3032102$, e-mail:caras@kocaeli.edu.tr;
carasgunduz@gmail.com
\\ $***$Gebze Institute of Technology, Department of Mathematics,\\
\c{C}ayirova Campus,  41400, Gebze Kocaeli,
 Turkey \\ Phone: $(+90262)6051389$ fax:$(+90262)6051388$, \\ e-mail:asonmez@gyte.edu.tr; ayse.sonmz@gmail.com   } 
\title{On lacunary statistically quasi-Cauchy sequences} 
\date{} 
\begin{document}
\newtheorem{Def}{Definition}
\newtheorem{Thm}{Theorem}
\newtheorem{Cor}{Corollary}
\newenvironment{pf}[1][Proof]{\noindent\textbf{#1.} }{\ \rule{0.5em}{0.5em}}
\newtheorem{Lem}{Lemma}
\maketitle
\begin{abstract}
The main object of this paper is to investigate lacunary
statistically ward continuity. We obtain some relations between
this kind of continuity and some other kinds of continuities. It
turns out that any lacunary statistically ward continuous real
valued function on a lacunary statistically ward compact subset
$E\subset{\textbf{R}}$ is uniformly continuous.
\end{abstract}
\section{Introduction}

A function $f:\textbf{R} \longrightarrow \textbf{R}$ is continuous
if and only if it preserves convergent sequences.  Using the idea
of continuity of a real function in terms of sequences, many kinds
of continuities were introduced by some researchers over the years,
and under different names many kinds of continuities have been
defined and investigated, not all but some of them we recall in
the following: slowly oscillating continuity
(\cite{CakalliSlowlyoscillatingcontinuity}), quasi-slowly
oscillating continuity (\cite{DikandCanak}), $\Delta$-quasi-slowly
oscillating continuity (\cite{CakalliNewkindsofcontinuities}, and
\cite{CakalliOnDeltaquasislowlyoscillatingsequences}), ward
continuity (\cite{CakalliForwardcontinuity}), $\delta$-ward
continuity (\cite{CakalliDeltaquasiCauchysequences}), statistical
ward continuity,  (\cite{CakalliStatisticalwardcontinuity}), and
$N_{\theta}$-ward continuity (\cite{CakalliNthetawardcontinuity})
which enabled some authors to obtain some characterizations of
uniform continuity in terms of sequences in the sense that a
function preserves either quasi-Cauchy sequences or slowly
oscillating sequences
(\cite{CakalliStatisticalquasiCauchysequences}, \cite{Vallin}, and
\cite{BurtonandColemanQuasiCauchySequences}).

The notion of $N_\theta$-convergence was introduced by  Freedman,
Sember and Raphael in
\cite{FreedmanandSemberandRaphaelSomecesarotypesummabilityspaces}.
Later on, the idea of lacunary statistically convergence was given
in \cite{FridyandOrhanLacunarystatisconvergence} based on the
notion of $N_\theta$-convergence. In
\cite{CakalliStatisticalwardcontinuity}, a subset $E$ of
\textbf{R} is called lacunary statistically ward compact if
whenever $(\alpha _{n})$ is a sequence of points in $E$ there is a
lacunary statistically-quasi-Cauchy subsequence of $(\alpha
_{n})$.

The aim of the present paper is to study the concept of lacunary
statistically ward continuity, and examine its properties.

\section{Preliminaries}
In this section, we remind basic constructions from the theory of
quasi-Cauchyness to lacunary quasi-Cauchyness.

Boldface letters $\boldsymbol{\alpha}$, $\bf{x}$, $\bf{y}$,
$\bf{z}$, ... will be used for sequences
$\boldsymbol{\alpha}=(\alpha_{n})$, $\textbf{x}=(x_{n})$,
$\textbf{y}=(y_{n})$, $\textbf{z}=(z_{n})$, ... of points in
$\textbf{R}$ for the sake of abbreviation. $s$ and $c$ will denote
the set of all sequences, and the set of convergent sequences of
points in \textbf{R}.

A sequence $(\alpha _{n})$ of points in $\textbf{R}$ is
quasi-Cauchy if $(\Delta \alpha _{n})$ is a null sequence, where
$\Delta \alpha _{n}=\alpha _{n+1}-\alpha _{n}$. These sequences
were named as quasi-Cauchy by Burton and Coleman \cite[page
328]{BurtonandColemanQuasiCauchySequences}, while they  were
called as forward convergent to $0$ sequences in
\cite{CakalliForwardcompactness} (see also \cite[page
226]{CakalliForwardcontinuity}). Quasi-Cauchy sequences arise in
diverse situations, and it is often difficult to determine whether
or not they converge, and if so, to which limit. It is easy to
construct a zero-one sequence such that the quasi-Cauchy average
sequence does not converge. The usual constructions have a
somewhat artificial feeling. Nevertheless, there are sequences
which seem natural, have the quasi-Cauchy property, and do not
converge. On the other hand, the sequence of averages of $0$ s and
$1$ s is always a quasi-Cauchy sequence: let $\textbf{x}:=
(x_{n})$ be a sequence such that for each nonnegative integer $n$,
$x_{n}$ is either $0$ or $1$. For each positive integer $n$ set
$a_{n}=\frac{x_{1}+x_{2}+...+x_{n}}{n}$ . Then $a_{n}$ is the
arithmetic mean average of the sequence up to time or position
$n$. Clearly for each $n$, $0 \leq a_{n} \leq 1$, and
$|a_{n+1}-a_{n}|  \leq {\frac{1}{n}}$. Thus $(a_{n})$ is a
quasi-Cauchy sequence.

The concept of statistical convergence is a generalization of the
usual notion of convergence that, for real-valued sequences,
parallels the usual theory of convergence (see \cite{Fridy},
\cite{MaioKocinac}, \cite{CakalliAstudyonstatisticalconvergence},
\cite{CasertaMaioKocinacStatisticalConvergenceinFunctionSpaces},
and \cite{CasertaandKocinacOnstatisticalexhaustiveness}). A
sequence $(x_{k})$ of points in $\textbf{R}$ is called
statistically convergent to an element $L$ of $\textbf{R}$ if
\[
\lim_{n\rightarrow\infty}\frac{1}{n}|\{k\leq n : |x_{k}-L|
\geq\varepsilon\}|=0,
\]
for every positive real number $\varepsilon$. In this case we
write $st-lim x_{k}=L$. The set of statistically convergent
sequences of points in $\textbf{R}$ is denoted by $S$. A sequence $(\alpha_{k})$ of points in $\textbf{R}$ is called to
be statistically quasi-Cauchy if $st-lim \Delta\alpha_{k}=0$. The set of statistically quasi-Cauchy sequences will be denoted by $\Delta S^{0}$.

A sequence $(x_{k})$ of points in $\textbf{R}$ is called lacunary
statistically convergent to an element $L$ of $\textbf{R}$ if
\[
\lim_{r\rightarrow\infty}\frac{1}{h_{r}}|\{k\in I_{r}: |x_{k}-L|
\geq\varepsilon\}|=0,
\]
for every positive real number $\varepsilon$ where
$I_{r}=(k_{r-1},k_{r}]$ and $k_{0}=0$,
$h_{r}:k_{r}-k_{r-1}\rightarrow \infty$ as $r\rightarrow\infty$
and $\theta=(k_{r})$ is an increasing sequence of positive
integers (\cite{FridyandOrhanLacunarystatisconvergence}). In this
case we write $S_{\theta}-lim x_{k}=L$. The set of lacunary
statistically convergent sequences of points in $\textbf{R}$ is
denoted by $S_{\theta}$. A sequence $(\alpha_{k})$ of points in $\textbf{R}$ is called to
be lacunary statistically quasi-Cauchy if $S_{\theta}-lim \Delta\alpha_{k}=0$. The set of lacunary statistically quasi-Cauchy sequences will be denoted by  $\Delta S^{0}_{\theta}$.

In the next section after reviewing lacunary statistically
sequentially continuity, we give the main results of our paper.

\section{Lacunary statistically ward continuity}
Connor and Grosse-Erdman (\cite{ConnorandGrosse}) gave sequential
definitions of continuity for real functions calling
$G$-continuity instead of $A$-continuity and their results covers
the earlier works related to $A$-continuity where a method of
sequential convergence, or briefly a method, is a linear function
$G$ defined on a linear subspace of $s$, denoted by $c_{G}$, into
$\textbf{R}$. A sequence $\textbf{x}=(x_{n})$ is said to be
$G$-convergent to $L$ if $\textbf{x}\in c_{G}$ and
$G(\textbf{x})=L$. A method $G$ is called regular if every
convergent sequence $\textbf{x}=(x_{n})$ is $G$-convergent with
$G(\textbf{x})=\lim \textbf{x}$. A method $G$ is called
subsequential if whenever $\textbf{x}$ is $G$-convergent with
$G(\textbf{x})=L$, then there is a subsequence $(x_{n_{k}})$ of
$\textbf{x}$ with $\lim_{k} x_{n_{k}}=L$. In particular, $\lim$
denotes the limit function $\lim \textbf{x}=\lim_{n}x_{n}$ on the
linear space $c$, $st-\lim$ denotes the statistical limit function
$st-\lim \textbf{x}=st-\lim_{n}x_{n}$ on the linear space $S$,
$S_{\theta}-\lim$ denotes the lacunary statistical limit function
$S_{\theta}-\lim \textbf{x}=S_{\theta}-\lim_{n}x_{n}$ on the
linear space $S_{\theta}$. Statistically sequential method is
regular without any restriction, however lacunary statistical
method is regular under the assumption that $\lim inf_{r}\;
q_{r}>1$. Statistical sequential method is subsequential, so is
lacunary statistical sequential method.

A subset $E$ of $\textbf{R}$ is called $G$-sequentially compact if
any sequence $\textbf{x}$ of points in $E$ has a $G$-sequentially
convergent subsequence $\textbf{z}$ such that $G(\textbf{z})\in
{E}$. A subset $E$ of $\textbf{R}$ is lacunary statistically
sequentially compact if any sequence $\textbf{x}$ of points in $E$
has a lacunary statistically convergent subsequence whose lacunary
statistical limit is in $E$. We see that this is a special case of
$G$-sequential compactness where $G=S_{\theta}-lim$. Lacunary
statistically sequentially compactness of a subset $E$ of
$\textbf{R}$ coincides not only with ordinary (sequential)
compactness, but also statistically sequentially compactness, and
lacunary statistically sequentially compactness.
 A function $f$ is called $G$-continuous at a point $u$ provided that whenever a sequence $\textbf{x}=(x_{n})$ of terms in the domain of $f$ is $G$-convergent to $u$, then the sequence $f(\textbf{x})=(f(x_{n}))$ is $G$-convergent to $f(u)$. Writing $G=S_{\theta}$, we get $S_{\theta}$-sequential continuity or lacunary statistically sequential continuity, explicitly we say that a real valued function $f$ defined on a subset $E$ of $R$ is called lacunary statistically sequentially continuous at a point $\alpha_{0}$ if it preserves lacunary statistically convergent sequences, i.e. $(f(\alpha_{k}))$ is a lacunary statistically convergent to $f(\alpha_{0})$ whenever $(\alpha_{k})$ is lacunary statistically convergent to $\alpha_{0}$.
Lacunary statistically sequentially continuity of a real valued
function defined on a subset of $\textbf{R}$ coincides with not
only ordinary (sequential) continuity, but also each one of the
continuities, statistically sequential continuity, and lacunary
statistically sequential continuity.

\textbf{Example 1.} Limit of the sequence of the ratios of
Fibonacci numbers converge to the golden mean. Lacunary sequential
method obtained via the sequence of Fibonacci numbers is a regular
method, i.e. $\theta=(k_{r})$ is the lacunary sequence defined by
writing $k_{0}=0$ and $k_{r}=F_{r+2}$ where $(F_{r})$ is the
Fibonacci sequence, i.e. $F_{1}=1$, $F_{2}=1$, $F_{r}= F_{r-1} +
F_{r-2}$ for $r\geq 3$ (see \cite{CakalliNthetawardcontinuity}). For this lacunary sequence $\theta=(k_{r})$, a real valued function defined on a subset of $\textbf{R}$ is $S_{\theta}$-sequentially (lacunary statistically sequentially) continuous if and only if it is ordinary (sequentially) continuous.

Using a similar idea to that of \cite{FridyandOrhanLacunarystatisconvergence} (see also \cite{CakalliLacunarystatisticalconvergenceintopgroups})  one can  easily find out the following inclusion properties between $\Delta S^{0}$ and $\Delta S^{0}_{\theta}$:\\
 (i) $\Delta S^{0}_{\theta} \subset {\Delta S^{0}}$ if and only if $lim\; sup \;q_{r} < \infty$ for any lacunary sequence $\theta$.\\
  (ii) $\Delta S^{0} \subset {\Delta S^{0}_{\theta}}$ if and only if $lim\; inf \;q_{r} > 1$ for any lacunary sequence $\theta$.

Combining these facts, for any lacunary sequence $\theta$, we have
that

(iii) $\Delta S^{0}_{\theta} = \Delta S^{0} $ if and only if
$1<lim\; inf \;q_{r}\leq lim\; sup \;q_{r} < \infty$;

(iv) $\Delta S^{0}_{\theta} = \Delta S^{0}$ if and only if
$S=S_{\theta}$;

(v) $\Delta S^{0}_{\theta} = \Delta S^{0}$ if and only if $|\Delta
N^{0}_{\theta}| = |\Delta \sigma^{0}_{1}|$;

(vi) $\Delta S^{0}_{\theta} = \Delta S^{0}$ if and only if $|\sigma_{1}|=N^{\theta}$\\
where $|\sigma_{1}|$ and $N^{\theta}$ denote the set of strongly
Cesaro convergent sequences and $N^{\theta}$-convergent sequences
of points in $\textbf{R}$, respectively. In the sequel, we will
always assume that $\lim inf_{r}\; q_{r}>1$.

In \cite{CakalliSequentialdefinitionsofconnectedness}, a non-empty
subset $E$ of $\textbf{R}$ is called {\em $G$-sequentially
connected} if there are no non-empty and disjoint $G$-sequentially
closed subsets $U$ and $V$ such that $A\subseteq U\bigcup V$, and
$A\cap U$ and $A\cap V$ are non-empty (see also
\cite{CakalliandMucukOnconnectednessviaasequentialmethod}). As far as $G$-sequentially connectedness is considered, we see
that  lacunary sequentially continuous image of any lacunary
sequentially connected subset of $\textbf{R}$ is lacunary
sequentially connected, so lacunary sequentially continuous image
of any interval is an interval. Furthermore taking $G=S_{\theta}$
it can be easily seen that a subset of $\textbf{R}$ is lacunary
sequentially connected if and only if it is connected in the
ordinary sense, and so is an interval.

Now we give the following interesting
examples which show emphasis the interest in different research
areas.

\textbf{Example 2.} Let $n$
be a positive integer. In a group of $n$ people, each person
selects at random and simultaneously another person of the group.
All of the selected persons are then removed from the group,
leaving a random number $n_{1} < n$ of people which form a new
group. The new group then repeats independently the selection and
removal thus described, leaving $n_{2} < n_{1}$ persons, and so
forth until either one person remains, or no persons remain.
Denote by $p_n$ the probability that, at the end of this iteration
initiated with a group of $n$ persons, one person remains. Then
the sequence $\textbf{p} = (p_{1}, p_{2}, · · ·, p_{n},...)$  is a
lacunary statistically quasi-Cauchy sequence, and $lim p_n$ does
not exist (see also \cite{WinklerMathematicalPuzzles}).

\textbf{Example 3.} Let $n$
be a positive integer. In a group of $n$ people, each person
selects independently and at random one of three subgroups to
which to belong, resulting in three groups with random numbers
$n_{1}$, $n_{2}$, $n_{3}$ of members; $n_{1} + n_{2} + n_{3} = n$.
Each of the subgroups is then partitioned independently in the
same manner to form three sub subgroups, and so forth. Subgroups
having no members or having only one member are removed from the
process. Denote by $t_{n}$ the expected value of the number of
iterations up to complete removal, starting initially with a group
of $n$ people. Then the sequence $(t_{1}, \frac{t_{2}}{2},
\frac{t_{3}}{3},...,\frac{t_{k}}{k},...)$ is a bounded
nonconvergent lacunary statistically quasi-Cauchy sequence (see also \cite{KeaneUnderstandingErgodicity}).

Now we state the definition of lacunary statistically ward
continuity in the following.
\begin{Def}
A real valued function $f$ defined on a subset $E$ of $\textbf{R}$
is called lacunary statistically ward continuous on $E$ if it
preserves lacunary statistically quasi-Cauchy sequences of points
in $E$, i.e. $(f(\alpha_{k}))$ is a lacunary statistically
quasi-Cauchy sequence whenever $(\alpha_{k})$ is a lacunary
statistically quasi-Cauchy sequences of points in $E$.
\end{Def}

We note that this definition of continuity cannot be obtained by
any $A$-continuity, i.e., by any summability matrix $A$, even by
the summability matrix $A=(a_{nk})$ defined by
$$a_{nk}=\frac{1}{h_{r}}\; for\; k=n+1,\;\;\;and\;\;\; a_{nk}=-\frac{1}{h_{r}}\; for\; k=n\;\; and\; a_{nk}=0\; otherwise$$
However, for this special summability matrix $A$, if
$A$-continuity of $f$ at the point $0$ implies lacunary
statistically ward continuity of $f$, then $f(0)=0$; and if
lacunary statistically ward continuity of $f$ implies
$A$-continuity of $f$ at the point $0$, then $f(0)=0$.

Sum of two lacunary statistically ward continuous functions is
lacunary statistically ward continuous, but product of lacunary
statistically ward continuous functions need not be lacunary
statistically ward continuous.

We give the following theorem the proof of which also can be
obtained by considering \cite[Theorem
9]{CakalliandHazarikaIdealquasiCauchysequences}.

\begin{Thm} \label{lacunarystatisticalwardcontinuityimplieslacunarystatiscontinuity} If a real valued function is lacunary statistically ward continuous on a subset $E$ of $\textbf{R}$, then it is lacunary statistically sequentially continuous on $E$.
\end{Thm}
\begin{pf} Suppose that $f$ is a lacunary statistically ward continuous function on a subset $E$ of $\textbf{R}$. Let $(x_{n})$ be a lacunary statistically quasi-Cauchy sequence of points in $E$. Then the sequence $$(x_1, x_0, x_2, x_0, x_3, x_0,..., x_{n-1}, x_0, x_n, x_0,...)$$ is a lacunary statistically quasi-Cauchy sequence. Since $f$ is lacunary statistically ward continuous, the sequence $$(y_n)=(f(x_{1}),f(x_{0}),f(x_{2}),f(x_{0}),...,f(x_{n}),f(x_{0}),...)$$ is a lacunary statistically quasi-Cauchy sequence. Therefore $S_{\theta}-\lim_{n\rightarrow \infty} \Delta y_n=0$. Hence  $S_{\theta}-\lim_{n\rightarrow \infty} [f(x_{n})-f(x_{0})]=0$. It follows that the sequence $(f(x_{n}))$ is lacunary statistically convergent to $f(x_{0})$. This completes the proof of the theorem.
\end{pf}

\begin{Cor} \label{lacunarystatisticalwardcontinuityimpliesordinarycontinuity} If a real valued function is lacunary statistically ward continuous on a subset $E$ of $\textbf{R}$, then it is ordinary continuous on $E$.
\end{Cor}

\begin{pf}
The proof follows immediately from the preceding theorem so is omitted.
\end{pf}

Now we prove the following theorem.
\begin{Thm} \label{TheoremuniformlycontinuousfunctiononEsendsquasiCauchytolacunarystatisticallquasiCauchy}
If a real valued function $f$ is uniformly continuous on a subset
$E$ of $\textbf{R}$, then $(f(x_{n}))$ is lacunary statistically
quasi-Cauchy whenever $(x_{n})$ is a quasi-Cauchy sequence of
points in $E$.
\end{Thm}
\begin{pf} Let $f$ be uniformly continuous on $E$. Take any quasi-Cauchy sequence $(x_{n})$ of points in $E$. Let $\varepsilon$ be any positive real number. Since $f$ is uniformly continuous, there exists a $\delta>0$ such that $|f(x)-f(y)|<\varepsilon$ whenever $|x-y|<\delta$. As $(x_{n})$ is a quasi-Cauchy sequence, for this $\delta$ there exists an $n_{0}\in{\textbf{N}}$ such that $|x_{n+1}-x_{n}|<\delta$ for $n\geq n_{0}$. Therefore $|f(x_{n+1})-f(x_{n})|<\varepsilon$ for $n\geq n_{0}$, so the number of indices $k$ for which $|f(x_{n+1})-f(x_{n})|\geq\varepsilon$ is less than $n_{0}$. Hence

$\lim_{r\rightarrow\infty}\frac{1}{h_{r} }|\{k\in I_{r}:|f(x_{n+1})-f(x_{n})|
\geq\varepsilon\}|\leq \lim_{r\rightarrow\infty}\frac{n_{0}}{h_{r} }=0$.\\
This completes the proof of the theorem.
\end{pf}

On the other hand, any continuous function on a compact subset $E$
of $\textbf{R}$ is uniformly continuous on $E$. It is also true
for a regular subsequential method $G$ that any lacunary
statistically ward continuous function on a $G$-sequentially
compact subset $E$ of $\textbf{R}$ is also uniformly continuous on
$E$ (see \cite{CakalliSequentialdefinitionsofcompactness}).
Furthermore, for lacunary statistically ward continuous functions
defined on a lacunary statistically ward compact subset of
\textbf{R}, we have the following.

\begin{Thm} \label{Lacunarystatisticalwardcontinuousfunctiononlacunarystatiswardcompactsubstisunifotmlycontinuous} Any lacunary statistically ward continuous real valued function on a lacunary statistically ward compact subset of $\textbf{R}$ is uniformly continuous.

\end{Thm}
\begin{pf} Let $E$ be a lacunary statistically ward compact subset $E$ of
$\textbf{R}$ and let $f:E\longrightarrow$ $\textbf{R}$ be a  lacunary
statistically ward continuous function on $E$. Suppose that $f$ is not
uniformly continuous on $E$ so that there exists an  $\varepsilon_{0} > 0$
such that for any $\delta >0$, there are $x, y \in{E}$ with $|x-y|<\delta$ but
$|f(x)-f(y)| \geq \varepsilon_{0}$. For each positive integer $n$, there exist
$\alpha_{n}$ and $\beta_{n}$ such that $|\alpha _{n}-\beta_{n}|<\frac{1}{h_{n}}$,
 and $|f(\alpha _{n})-f(\beta_{n})|\geq \varepsilon_{0}$. Since $E$ is lacunary
 statistically ward compact, there exists a lacunary statistically quasi-Cauchy
 subsequence $(\alpha _{n_{k}})$ of the sequence $(\alpha _{n})$. It is clear that
  the corresponding subsequence $(\beta_{n_{k}})$ of the sequence $(\beta_{n})$ is
   also lacunary statistically quasi-Cauchy, since $(\beta_{n_{k+1}}-\beta_{n_{k}})$ is lacunary statistically convergent to $0$ which follows from the following lines: for each $\varepsilon$,
 $\{k\in I_{r}: |\beta_{n_{k+1}}-\beta_{n_{k}}| \geq \varepsilon\}
 \subset {\{k\in I_{r}: |\beta_{n_{k+1}}-\alpha_{n_{k+1}}|\geq
 \frac{\varepsilon}{3}\} }$ $$\; \; \;\cup \{k\in I_{r}: |\alpha_{n_{k+1}}-\alpha_{n_{k}}|
 \geq  \frac{\varepsilon}{3} \} \cup \{k\in I_{r}: |\alpha_{n_{k}}-\beta_{n_{k}}| \geq \frac{\varepsilon}{3}\}.$$ It follows from this inclusion that \\
$|\{k\in I_{r}: |\beta_{n_{k+1}}-\beta_{n_{k}}| \geq
\varepsilon\}|$ $\leq {|\{k\in I_{r}:
|\beta_{n_{k+1}}-\alpha_{n_{k+1}}|\geq
\frac{\varepsilon}{3}\}}|$\\ $+ |\{k\in I_{r}:
|\alpha_{n_{k+1}}-\alpha_{n_{k}}| \geq \frac{\varepsilon}{3}\}|+$
$|\{k\in I_{r}: |\alpha_{n_{k}}-\beta_{n_{k}}| \geq
\frac{\varepsilon}{3}\}|$.

Hence

$\lim_{r\rightarrow\infty}\frac{1}{h_{r} }|\{k\in I_{r}:
|\beta_{n_{k+1}}-\beta_{n_{k}}| \geq \varepsilon\}|$ $\leq
{\lim_{r\rightarrow\infty}\frac{1}{h_{r} }|\{k\in I_{r}:
|\beta_{n_{k+1}}-\alpha_{n_{k+1}}|\geq \frac{\varepsilon}{3}
\}}|+$ $\lim_{r\rightarrow\infty}\frac{1}{h_{r} }|\{k\in I_{r}:
|\alpha_{n_{k+1}}-\alpha_{n_{k}}| \geq \frac{\varepsilon}{3}\}|+$
$\lim_{r\rightarrow\infty}\frac{1}{h_{r} }|\{k\in I_{r}:
|\alpha_{n_{k}}-\beta_{n_{k}}| \geq \frac{\varepsilon}{3}\}|=0$.

On the other hand, it follows from the equality
$$\alpha_{n_{k+1}}-\beta_{n_{k}}=\alpha_{n_{k+1}}-\alpha_{n_{k}}+\alpha_{n_{k}}-\beta_{n_{k}}$$
that the sequence $(\alpha_{n_{k+1}}-\beta_{n_{k}})$ is lacunary
statistically convergent to $0$. Hence the sequence
$$(a_{n_{1}}, \beta_{n_{1}}, \alpha _{n_{2}}, \beta_{n_{2}}, \alpha _{n_{3}}, \beta_{n_{3}},..., \alpha _{n_{k}}, \beta_{n_{k}},...)$$
is lacunary statistically quasi-Cauchy. But the transformed
sequence
$$(f(\alpha _{n_{1}}), f(\beta_{n_{1}}), f(\alpha _{n_{2}}), f(\beta_{n_{2}}), f(\alpha _{n_{3}}), f(\beta_{n_{3}}),..., f(\alpha _{n_{k}}), f(\beta_{n_{k}}),...)$$
is not lacunary statistically quasi-Cauchy. Thus $f$ does not
preserve lacunary statistically quasi-Cauchy sequences. This
contradiction completes the proof of the theorem.
\end{pf}

\begin{Cor} \label{Lacunarytatiswardcontinuousfunctiononboundedsubsetisunifotmlycontinuous} If a real valued function is lacunary statistically ward continuous on a bounded subset $E$ of \textbf{R}, then it is uniformly continuous on $E$.
\end{Cor}

\begin{pf}
The proof follows from the preceding theorem and Theorem 3 in
\cite[page 1622]{CakalliStatisticalquasiCauchysequences}.
\end{pf}

\begin{Cor} \label{LacunaryatiswardcontinuousfunctiononNthetawardcompactsubsetisunifotmlycontinuous} If a real valued function is lacunary statistically ward continuous on an $N_{\theta}$-ward compact subset $E$ of \textbf{R}, then it is uniformly continuous on $E$.
\end{Cor}

\begin{pf}
The proof follows from Theorem
\ref{Lacunarystatisticalwardcontinuousfunctiononlacunarystatiswardcompactsubstisunifotmlycontinuous}
and \cite[Theorem 3.3]{CakalliNthetawardcontinuity}.
\end{pf}

We give below that any real valued lacunary statistically-ward
continuous function defined on an interval is uniformly
continuous. First we give the following lemma.

\begin{Lem} \label{Lemmaanypairofsequencehasalamdastatisticallystatisticallyquasi}
If $(\xi_{n}, \eta_{n})$ is a sequence of ordered pairs of points
in an interval such that $\lim_{n\rightarrow\infty}
|\xi_{n}-\eta_{n}|=0$, then there exists a  lacunary statistically
quasi-Cauchy sequence $(\alpha_{n})$ with the property that for
any positive integer $i$ there exists a positive integer $j$ such
that $(\xi_{i}, \eta_{i})=(\alpha_{j-1}, \alpha_{j})$.
\end{Lem}

\begin{Thm} \label{TheoremLamdastatisticallywardcontinuityonanintervalimpliesuniformlycontinuity}
If a real valued function defined on an interval $E$ is lacunary
statistically ward continuous, then it is uniformly continuous.
\end{Thm}
\begin{pf}
Suppose that $f$ is not uniformly continuous on $E$. Then there is
an $\varepsilon_{0} > 0$ such that for any $\delta > 0$ there
exist $x, y \in {E}$ with $|x - y| < \delta$ but $| f (x) - f
(y)|\geq {\varepsilon_{0}}$. For every $n \in{\textbf{N}}$ fix
$\xi_{n}$, $\eta_{n}\in{E}$ with $|\xi_{n} - \eta_{n}| <
\frac{1}{n}$ and $| f (\xi_{n}) - f (\eta_{n})| \geq
\varepsilon_{0}$. By Lemma
\ref{Lemmaanypairofsequencehasalamdastatisticallystatisticallyquasi},
there exists a lacunary  statistically quasi-Cauchy sequence
$(\alpha_{i})$ such that for any integer $i \geq {1}$ there exists
a $j$ with $\xi_{i} = \alpha_{j}$ and $\eta_{i} = \alpha_{j+1}$.
This implies that $| f (\alpha_{j+1} ) - f (\alpha_{j})| \geq
\varepsilon_{0}$; hence $( f (\alpha_{i} ))$ is not  lacunary
statistically quasi-Cauchy. Thus $f$ does not preserve lacunary
statistically quasi-Cauchy sequences. This completes the proof of
the theorem.
\end{pf}

Since the sequence constructed in Lemma
\ref{Lemmaanypairofsequencehasalamdastatisticallystatisticallyquasi}
is also quasi-Cauchy, we see that the statement $(f(x_{n}))$ is
lacunary statistically quasi-Cauchy whenever $(x_{n})$ is
quasi-Cauchy sequence of points in $E$ implies the uniform
continuity of $f$ on $E$. Now combining Theorem
\ref{TheoremuniformlycontinuousfunctiononEsendsquasiCauchytolacunarystatisticallquasiCauchy}
with this observation we have the following result.

\begin{Cor}
Let $f$ be a real valued function defined on an interval $E$. Then
$f$ is uniformly continuous on $E$ if and only if $(f(x_{n}))$ is
lacunary statistically quasi-Cauchy whenever $(x_{n})$ is
quasi-Cauchy sequence of points in $E$.
\end{Cor}

\begin{Cor}
Let $f$ be a real valued function defined on an interval $E$. Then the following statements are equivalent:\\
(a)   if $(f(x_{n}))$ is lacunary statistically quasi-Cauchy whenever $(x_{n})$ is quasi-Cauchy sequence of points in $E$.\\
(b)   if $(f(x_{n}))$ is $N_{\theta}$ quasi-Cauchy whenever
$(x_{n})$ is quasi-Cauchy sequence of points in $E$.
\end{Cor}

\begin{pf}
The proof follows from Theorem
\ref{TheoremLamdastatisticallywardcontinuityonanintervalimpliesuniformlycontinuity},
and \cite[Theorem 1 and Theorem
2]{CakalliandKaplanAstudyonNthetawardcontinuity} so is omitted.
\end{pf}

\begin{Cor}
If a real valued function defined on an interval is lacunary
statistically ward continuous, then it is ward continuous.
\end{Cor}

\begin{pf}
The proof follows from Theorem
\ref{TheoremLamdastatisticallywardcontinuityonanintervalimpliesuniformlycontinuity},
and \cite[Theorem 5]{CakalliStatisticalquasiCauchysequences} so it
is omitted.
\end{pf}
\begin{Cor}
If a real valued function defined on an interval is lacunary
statistically ward continuous, then it is slowly oscillating
continuous.
\end{Cor}
\begin{pf}
The proof follows from Theorem
\ref{TheoremLamdastatisticallywardcontinuityonanintervalimpliesuniformlycontinuity},
and \cite[Theorem 5]{CakalliStatisticalquasiCauchysequences} so it
is omitted.
\end{pf}

\section{Conclusion}
The concept of continuity and any concept involving continuity
play a very important role not only in pure mathematics but also
in other branches of sciences involving mathematics especially in
computer science, information theory, biological science. In this
paper, the concept of lacunary statistically ward continuity of a
real function is investigated. In this investigation we have
obtained theorems related to lacunary statistically ward
continuity, $N_{\theta}$-ward continuity, slowly oscillating continuity, ward continuity, uniform continuity, lacunary
statistically ward compactness, $N_{\theta}$-ward compactness, boundedness, and compactness. Although we have not given in the paper, it is not difficult to find out that uniform limit of a sequence of lacunary statistically ward continuous
functions is lacunary statistically ward continuous, and
that the set of all lacunary statistically ward continuous
functions on a subset $E$ of $\textbf{R}$ is a closed subset of
the set of all continuous functions on $E$.

For a further study, we suggest to investigate lacunary
statistically quasi-Cauchy sequences of fuzzy points, and lacunary
statistically ward continuity for the fuzzy functions (see
\cite{SavasRemarkondoublelacunarystatisticalconvergenceoffuzzynumbers},
\cite{SavasOnasymptoticallylacunarystatisticalequivalentsequencesoffuzzynumbers},
\cite{MursaleenandMohiuddineOnlacunarystatisticalconvergencewithrespecttotheintuitionisticfuzzynormedspace},
and \cite{CakalliandPratul} for the definitions and  related
concepts in fuzzy setting). However due to the change in settings,
the definitions and methods of proofs will not always be analogous
to those of the present work. For another further study we suggest
to investigate lacunary statistically quasi-Cauchy sequences of
double sequences, and lacunary statistically ward double
continuity to find out whether lacunary statistically ward double
continuity coincides with lacunary  statistically ward (single)
continuity or not (see
\cite{PattersonandSavasLacunarystatisticalconvergenceofdoublesequences} for the definitions and related concepts in the double case).

\end{document}